\documentclass[a4paper,11pt]{amsart}
\usepackage[all]{xy}        % Commutative Diagrams
\CompileMatrices            % Faster Tex-ing
\UseTips                    % Use Computer Modern Arrowheads
\usepackage[german,english]{babel}   %Foreign language support
\usepackage[bookmarks=true]{hyperref}   % Hyperref in DVI and PDF (like HTML Links between sections)
\usepackage{amssymb}
\usepackage{xspace}
\theoremstyle{plain}
\newtheorem{theorem}{Theorem}[section]
\newtheorem*{theorem*}{Theorem}
\newtheorem{prop}[theorem]{Proposition}
\newtheorem{cor}[theorem]{Corollary}
\newtheorem{lem}[theorem]{Lemma}

\theoremstyle{definition}
\newtheorem{defn}[theorem]{Definition}

\theoremstyle{remark}
\newtheorem{notation}[theorem]{Notation}
\newtheorem{remk}[theorem]{Remark}

\newtheorem{examples}[theorem]{Examples}

\numberwithin{theorem}{section}

\newcommand{\enm}[1]{\ensuremath{#1}}          % Shortcuts
\newcommand{\op}[1]{\operatorname{#1}}
\newcommand{\cal}[1]{\mathcal{#1}}
\newcommand{\wt}[1]{\widetilde{#1}}
\newcommand{\wh}[1]{\widehat{#1}}
\renewcommand{\bar}[1]{\overline{#1}}

\newcommand{\CC}{\enm{\mathbb{C}}}             % All Number domains easily accssable

\newcommand{\FF}{\enm{\mathbb{F}}}

\newcommand{\PP}{\enm{\mathbb{P}}}

\newcommand{\Dd}{\enm{\cal{D}}}

\newcommand{\Hh}{\enm{\cal{H}}}
\newcommand{\Ii}{\enm{\cal{I}}}

\newcommand{\Ll}{\enm{\cal{L}}}
\newcommand{\Mm}{\enm{\cal{M}}}
\newcommand{\Nn}{\enm{\cal{N}}}
\newcommand{\Oo}{\enm{\cal{O}}}

\renewcommand{\phi}{\varphi}        % Dont know how to not loose the original ones???
\renewcommand{\theta}{\vartheta}
\renewcommand{\epsilon}{\varepsilon}

\newcommand{\Ann}{\op{Ann}}         % Standard Operators

\newcommand{\Spec}{\op{Spec}}

\newcommand{\Hom}{\op{Hom}}
\newcommand{\Ext}{\op{Ext}}

\newcommand{\id}{\op{id}}
\newcommand{\dirlim}{\varinjlim}
\newcommand{\invlim}{\varprojlim}

\newcommand{\tensor}{\otimes}         % Symbols with meaning

\newcommand{\xd}[1][x]{\enm{#1_1,\ldots,#1_d}}

\renewcommand{\to}[1][]{\xrightarrow{\ #1\ }}

\newcommand{\into}{\hookrightarrow}

\newcommand{\usc}[1][m]{\underline{\phantom{#1}}}

\newcommand{\defeq}{\stackrel{\scriptscriptstyle \op{def}}{=}}

\newcommand{\ie}{\textit{i.e.}\ }           % i.e. in italics and with proper spacing afterwards
\newcommand{\cf}{\textit{cf.}\ }

\newcommand{\F}[1][]{F^{#1*}}

\newcommand{\Czech}{\v{C}ech\xspace}               % The Czech guy
\newcommand{\comment}[1]{}

\newcommand{\Fful}[1]{F^{\infty}{#1}}
\newcommand{\Fnil}[1]{{#1}_{\op{nil}}}         % F-nilpotent submodule
\newcommand{\Fred}[1]{{#1}_{\op{red}}}         % F-reduced quotient
\newcommand{\Ffred}[1]{{#1}_{\op{fred}}}       % F-full reduced subquotient

\begin{document}
\title{The intersection homology $D$--module in finite characteristic.}
\author{Manuel Blickle}
\address{Universit\"at Essen, FB6 Mathematik, 45117 Essen,
Germany} \email{manuel.blickle@uni-essen} \keywords{D-modules,
unit F-modules, Intersection Homology, Tight Closure}
\subjclass[2000]{14B15,13N10}

\maketitle

\begin{abstract}
For $Y$ a closed normal subvariety of codimension $c$ of a smooth
$\CC$--variety $X$, Brylinski and Kashiwara showed that the local cohomology
module $\Hh^c_Y(X,\Oo_X)$ contains a unique simple $\Dd_X$--submodule, denoted
by $\Ll(Y,X)$. In this paper the analogous result is shown for $X$ and $Y$
defined over a perfect field of finite characteristic. Moreover, a local
construction of $\Ll(Y,X)$ is given, relating it to the theory of tight
closure. From the construction one obtains a criterion for the
$\Dd_X$--simplicity of $\Hh^c_Y(X,\Oo_X)$.
\end{abstract}

%\tableofcontents

\section{Introduction}
Let $Y$ be a closed codimension $c$ subvariety of the smooth $\CC$
variety $X$ and let $Z$ be the singular locus of $Y$. Denote by
$\Dd_X$ the sheaf or differential operators on $X$. In
\cite[Proposition 8.5]{Bry.Kash}, Brylinski and Kashiwara show the
existence (and usefulness) of a unique holonomic $\Dd_X$ module
$\Ll=\Ll(Y,X)$ satisfying the properties
\begin{eqnarray*}
    \Ll|_{X-Z} \cong \Hh^c_{Y-Z}(X-Z,\Oo_{X-Z}) \\
    \Hh_Z^0(\Ll) = \Hh_Z^0(\Ll^*) = 0,\qquad
\end{eqnarray*}
where the star stands for duality of holonomic $\Dd$--modules and $\Hh^i_Y$
denotes the higher derived sections with support in $Y$. The proof of this
result is rather formal and uses duality theory for holonomic $\Dd_X$--modules.
Furthermore, they show that $\Ll(Y,X)$ is the unique simple, selfdual holonomic
$\Dd_X$--module agreeing with $\Hh_Y^c(X,\Oo_X)$ on $X-Z$. This result is
obtained by showing that $\Ll(Y,X)$ corresponds, via the Riemann--Hilbert
correspondence, to the intersection homology complex $\pi_Y$ of middle
perversity, which, by construction, is simple and selfdual. All these
constructions, such as holonomicity, duality and the Riemann--Hilbert
correspondence completely rely on characteristic zero -- on analytic techniques
even, if one is strict.

The question answered in this paper is: What is the situation if $X$ is defined
over a field of positive characteristic? Somewhat surprisingly, the existence
of a unique simple $\Dd_X$--submodule $\Ll(Y,X)$ can be proved almost
independent of the characteristic. The key ingredient -- the proof of which is
characteristic dependent though -- is that $\Hh^c_Y(X,\Oo_X)$ has finite length
as a $\Dd_X$--module. This is guaranteed by holonomicity in characteristic $0$
and by \cite[Theorem 5.7]{Lyub} in positive characteristic, respectively.

We state the result and sketch the simple argument -- for a
complete proof refer to Theorem \ref{thm.main}.
\begin{theorem*}
    Let $X$ be a smooth $k$--variety and let $Y$ be a closed
    irreducible subvariety of codimension $c$. Then $\Hh^c_Y(X,\Oo_X)$ has
    a unique simple $\Dd_X$--submodule $\Ll(Y,X)$. Furthermore, $\Ll(Y,X)$ agrees
    with $\Hh^c_Y(X,\Oo_X)$ on $X- \op{Sing} Y$.
\end{theorem*}
\begin{proof}(Sketch)
Since $H^c_Y(X,\Oo_X)$ has finite length as a $\Dd_X$--module it has some
simple non-zero $\Dd_X$--submodule $\Ll$. Denote the inclusion $X'\defeq
X-\op{Sing Y} \subseteq X$ by $i$ and write $Y'$ for $Y-\op{Sing} Y$. One sees
easily that the restriction of $\Ll$ to $X'$ is nonzero. As the restriction of
$\Hh^c_Y(X,\Oo_X)$ is equal to $\Hh^c_{Y'}(X',\Oo_{X'})$, and by smoothness of
$Y'$, the latter is $\Dd_X$--simple it follows that $\Ll|_{X'} =
\Hh^c_Y(X,\Oo_X)|_{X'}$. Since this holds for any simple submodule of
$\Hh^c_Y(X,\Oo_X)$ the same argument shows that any two such have nonzero
intersection, thus they are equal. This shows the uniqueness.
\end{proof}
This existence proof gives very little information about the
concrete structure of $\Ll(Y,X)$. Even in characteristic zero, to
explicitly determine $\Ll(Y,X)$ is difficult. The best results in
this case are due to Vilonen \cite{Vil} for $Y$ a complete
intersection with an isolated singularity. He uses analytic
techniques to characterize the sections of $\Hh^c_Y(X,\Oo_Y)$
belonging to $\Ll(Y,X)$. They are precisely the ones vanishing
under a certain residue map. Furthermore he gives a canonical
generator, the canonical class associated to $Y \subseteq X$, for
$\Ll(Y,X)$ in this case.

To explicitly determine $\Ll(Y,X)$ in positive characteristic is
the main purpose of this paper. The strategy is to use the
Frobenius instead of the differential structure. This substitution
is justified by the close relationship of so called unit
$\Oo_X[F^e]$--structures and $\Dd_X$--structures, described in
\cite{Lyub,Manuel.DRF,Em.Kis2}. Our construction is local in
nature. If we denote by $R$ and $A=R/I$ the local rings of $X$ and
$Y$ at a point $x \in Y$, we roughly show the following, for
precise statements see Section \ref{sec.LocConstLl}.
\begin{theorem*}
Let $R$ be regular, local and $F$--finite. Let $A=R/I$ be a normal
domain. Then the unique simple $D_R$--submodule, $\Ll(A,R)$, of
$H^c_I(R)$ is dual to the unique simple $A[F^e]$--module quotient
of $H^d_m(A)$.
\end{theorem*}
The duality we are referring to is an extension of Matlis duality
incorporating Frobenius actions. Furthermore, the construction is
explicit enough to identify (non canonical) generators for
$\Ll(A,R)$. What we have gained is that the unique simple
$A[F^e]$--module quotient of $H^d_m(A)$ is well studied and fairly
well understood; it is the quotient of $H^d_m(A)$ by the tight
closure of zero, $0^*_{H^d_m(A)}$. The vanishing of
$0^*_{H^d_m(A)}$ is governed by $F$--rationality of $A$, which is
a positive characteristic analog of rational singularities. As a
consequence of this connection we obtain the following
$D_R$--simplicity criterion for $H^c_I(R)$:
\begin{theorem*}
Let $R$ be regular, local and $F$--finite. Let $A=R/I$ be a
Cohen--Macaulay domain of codimension $c$. Then, if $A$ is
$F$--rational then $H^c_I(R)$ is $D_R$--simple.
\end{theorem*}
More precise simplicity criteria for $H^c_I(R)$ are given in
Section \ref{sec.SimCritH}.

The paper is structured as follows. In Sections 2 and 3 we recall and further
develop some necessary machinery from the theory of $R[F^e]$--modules and tight
closure. As the techniques used later are local in nature, the notation
reflects this and we mainly speak of rings and ideals instead of schemes and
their sub-schemes. In these sections we do not concretely deal with the
applications to constructing $\Ll(Y,X)$ but derive general results which
constitute the technical underpinning of what follows. As a notable byproduct
we answer a question posed by Lyubeznik showing that minimal roots exist for
finitely generated unit $R[F^e]$--modules for any regular, local ring $R$. In
\cite{Lyub} this was only shown in the complete case.

Section 4 contains the main results discussed above and generalizations
thereof. Furthermore, as an application to tight closures theory we show that
the parameter test module commutes with localization. We finish this section
with a complete characterization of $D_R$--simplicity in the case of curves,
providing a finite characteristic analog of results of Yekutieli
\cite{Yeku.ResDif} and S.\ P.\ Smith \cite{SmithSP.IntHom}.

Finally we remark that the substitution of the $\Dd_X$--module
structure by a unit $\Oo_X[F^e]$--structure in the study of
$\Ll(A,R)$ enables one to place $\Ll(Y,X)$, in finite
characteristic, in the context of a Riemann--Hilbert type
correspondence. That such a correspondence exists is recent work
of Emerton and Kisin \cite{Em.Kis,Em.Kis2}, where an equivalence
(on the level of derived categories) of the category of finitely
generate unit $\Oo_X[F^e]$--modules and the category of
constructible $\FF_{p^e}$--sheaves is developed. Within this
correspondence, the simple unit $\Oo_X[F^e]$--module $\Ll(Y,X)$
constructed here does indeed correspond to certain middle
extensions on the constructible $\FF_p$--site. These connections
will not be discussed here but should appear in the final version
of \cite{Em.Kis2}, and are outlined in \cite{Em.Kis3}.

It is a pleasure to thank my advisor Karen Smith for her expertise and
encouragement during my dissertation \cite{Manuel.PhD}, which contains most of
the results presented here. Furthermore, I thank Matt Emerton, Brian Conrad and
Gennady Lyubeznik for valuable conversations on various parts of this paper.

\section{Background on $R[F^e]$--modules}
Throughout this paper, $R$ denotes a noetherian ring of dimension
$n$ containing a field $k$ of positive characteristic $p$, unless
stated otherwise. For an ideal $I$ of height $c$ we denote the
quotient $R/I$ by $A$. This is a ring of dimension $d=n-c$. In
general we assume that $R$ is regular and $F$--finite, \ie $R$ is
a finite module over its subring of $p$th powers.

The \emph{(absolute) Frobenius map} on $R$, \ie the ring map sending each
element to its $p$th power, is denoted by $F=F_R$. The associated map on
$X=\Spec R$ is denoted by the same letter $F=F_X$.

If $\Mm$ is an $R$--module, then $\Mm^e$ denotes the $R$--$R$--bimodule, which,
as a left module is just $\Mm$, but with right structure twisted by the $e$th
iterate of the Frobenius, \ie for $r \in R$ and $m \in \Mm$ one has $m \cdot r
= r^{p^e}m$. With this notation Psekine and Szpiro's \emph{Frobenius functor}
is defined as $F^{*}(\Mm)=R^1 \tensor \Mm$. Clearly, $F^{*}$ commutes with
direct limits and direct sums. If $R$ is regular, $F^{*}$ is flat; therefore it
commutes with finite intersections. The flatness of $F^*$ in the regular case
is where the theory draws its power from. The same is valid for higher powers
of the Frobenius and clearly we have $(F^e)^*=(F^*)^e$ which we denote by
$\F[e]$.

We review the definition and basic properties of modules with Frobenius action.
Since we are being extremely brief with this here, we advise the reader with no
prior exposure to first consult Section 2 of \cite{Manuel.DRF}. For a thorough
introduction see for example \cite{Manuel.PhD}, Chapter 2, or
\cite{Em.Kis,Lyub}.

\begin{defn}
An \emph{$R[F^e]$--module} is an $R$--module $\Mm$ together with an
$R$--linear map
\[
\theta^e: \F[e]\Mm = R^e \tensor \Mm \to \Mm.
\]
If $\theta^e$ is an isomorphism, then $(\Mm,\theta^e)$ is called a
\emph{unit $R[F^e]$--module}.
\end{defn}
By adjointness these maps $\theta^e \in \Hom(\F[e]\Mm,\Mm)$ are in one-to-one
correspondence with maps $F^e_\Mm \in \Hom(\Mm,F_*^e\Mm)$ where $F^e_\Mm(m) =
\theta^e_\Mm(1 \tensor m)$. Therefore, an $R[F^e]$--module is nothing but a
module over the non-commutative ring
\[
    R[F^e]=\frac{R\langle F^e\rangle}{r^{p^e}F^e - F^{e}r},
\]
where $F^e$ acts on $M$ via $F^e_\Mm$. Then the category of $R[F^e]$--modules,
$R[F^e]$--mod, is the category of left modules over this ring $R[F^e]$. As the
module category over an associative ring, $R[F^e]$--mod is an abelian category.
The category of unit $R[F^e]$--modules, $uR[F^e]$--mod, is the full subcategory
whose objects are those $R[F^e]$--modules which are unit. Since $R$ is regular,
the resulting flatness of $\F[e]$ implies that this is also an abelian
category. If $\Nn$ is an $R$--submodule of the $R[F^e]$--module $(\Mm,
\theta^e, F^e)$, then we denote for convenience $\theta^e(\F[e]\Nn)$ by
$F^e(\Nn)$, which is just the $R$-submodule of $\Mm$ generated by all $F^e(n)$
for $n \in \Nn$.
\begin{defn}
    An $R[F^e]$--module $(\Mm, \theta^e)$ is called finitely
    generated if it is a finitely generated module over the ring
    $R[F^e]$.
\end{defn}

Let $\phi: M \to \F[e]M$ be an $R$--linear map. Consider the directed limit of
the system of Frobenius powers of this map
\[
\Mm={\dirlim (M \to[\phi] \F[e]M
\to[F^*\phi] \F[2e]M \to \cdots \quad)}
\]
which carries a natural unit $R[F^e]$--module structure. If a unit
$R[F^e]$--module $(\Mm,\theta^e)$ arises in such a fashion one calls $\phi$ a
\emph{generator} of $(\Mm,\theta^e)$. If $M$ is finitely generated it is called
a finite generator, and if, in addition, $\phi$ is injective, then $\Mm$ is
called a \emph{root} of $\Mm$. In this case one identifies $M$ with its
isomorphic image in $\Mm = \dirlim \F[er]M$. Thus, a root of a unit
$R[F]$--module $\Mm$ is a finitely generated $R$--submodule $M$, such that $M
\subseteq RF^e(M)$ and $\Mm = \bigcup_r RF^{er}(M) = R[F^e]M$. A key
observation is the following proposition, see \cite[Proposition
2.5]{Manuel.DRF} or \cite{Em.Kis} for proof:
\begin{prop}
    Let $R$ be regular. A unit $R[F^e]$--module $(\Mm,\theta^e)$ is finite\-ly
    generated  if and only if $\Mm$ has a root.
\end{prop}
With this at hand one can easily show that the category of finitely generated
unit $R[F^e]$--modules is an abelian subcategory with ACC of the category of
$R[F^e]$--modules which is closed under extensions. Significantly more work
(for the second part) is involved in showing the next important theorem, found
as Proposition 2.7 and Theorem 3.2 in \cite{Lyub}.
\begin{theorem}\label{thm.RFfinite}
    Let $R$ be regular and let $\Mm$ be a finitely generated unit $R[F^e]$--module.
    Then $\Mm$ has ACC in the category of unit $R[F^e]$--modules

    If $R$ is also a finitely generated algebra over a regular local ring,
    then $\Mm$ has DCC, \ie $\Mm$ has finite length as a unit
    $R[F^e]$--module.
\end{theorem}

\begin{examples}\label{ex.rfmods} Standard examples of unit $R[F^e]$--modules
are: (1) $R$ itself via the natural isomorphism $R^e\tensor_R R \cong R$. (2) A
localization $S^{-1}R$ of $R$ via the natural map $R^e \tensor S^{-1}R \to
S^{-1}R$ whose inverse is the map $rs^{-1} \mapsto s^{p^e-1}r \tensor s^{-1}$.
(3) The local cohomology modules, $H^i_I(R)$, of $R$ with support in an ideal
$I$ which obtain their unit structure via the \Czech complex which consists of
localizations of $R$.
\end{examples}
Note that any proper nonzero ideal $I \subseteq R$ is an $R[F^e]$--submodule
but not unit. Thus $R$ is a \emph{simple} unit $R[F^e]$--module.

\subsection{Base change} Let $\pi: R \to S$ be a map of rings. The base
change functor $\pi^* = S \tensor_R \usc $ extends to a functor from (unit)
$R[F^e]$--modules to (unit) $S[F^e]$--modules. For $(\Mm,\theta^e)$, the (unit)
$S[F^e]$--module structure on $S\tensor_R \Mm$ is given by
\[
    S^e \tensor_S S \tensor_R \Mm \cong S \tensor_R R^e \tensor_R
    \Mm \to[\id_S \tensor \theta^e] S \tensor_R \Mm.
\]
Clearly, this is an isomorphism if and only if $\theta^e$ is an isomorphism.
For easy reference we record some properties of base change, the easy proofs
are left to the reader.
\begin{prop}
    Let\/ $R \to S$ be a map of rings. Let\/ $\Mm$ be a finitely generated
    unit\/
    $R[F^e]$--module with generator\/ $M$. (1) $S \tensor M$ is a generator of
    the finitely generated unit\/ $S[F^e]$--module\/ $S \tensor\Mm$. (2) If $R$
    and $S$ are regular, then the image of\/ $S \tensor M$ in\/ $S \tensor
    \Mm$ is a root of $S \tensor \Mm$. (3) If\/ $R \to S$ is also flat and $M$
    is a root of $\Mm$, then\/ $S \tensor M$ itself is a root of\/ $S \tensor \Mm$.
    (4) If\/ $R \to S$ is faithfully flat, then a submodule\/ $M$ of\/ $\Mm$ is
    a root of\/ $\Mm$ if and only if\/ $S \tensor M$ is a root of\/ $S \tensor \Mm$.
\end{prop}

\subsection{Restriction} Still fixing the data of a map of rings $\pi: R \to S$, any
$S[F^e]$--module $(\Nn,\theta^e)$ naturally carries an $R[F^e]$--module
structure because $\pi$ induces a ring homomorphism $R[F^e] \to S[F^e]$. Note
that in general, a unit $S[F^e]$--module, viewed as an $R[F^e]$--module (really
$\pi_*\Nn$), is \emph{not} unit. What is the case is that restriction of
scalars preserves the unit property, if and only if the relative Frobenius
$F^e_{S/R}: R^e \tensor_R S \to S^e$ sending $r \tensor s$ to $\pi(r)s^{p^e}$
is an isomorphism of $R$--$S$--bimodules. The following proposition summarizes
some cases where this happens:
\begin{prop}\label{prop.loccompet}
Let $\pi : R \to S$ be a map of rings. In the following cases, the relative
Frobenius $F^e_{S/R}: R^e \tensor_R S \to S^e$ is an isomorphism:
\begin{enumerate}
    \item $S$ is the localization of $R$ at some multiplicative set $T \subseteq R$,
    and $\pi$ is the localization map.
    \item $R \to S$ is \'etale.
    \item $R$ is regular local and $F$--finite, and $S$ is the $I$--adic completion
    of $R$ with respect to some ideal $I$ of $R$.
\end{enumerate}
In all these cases it follows that a unit $S[F^e]$--module is also unit as a
$R[F^e]$--module.
\end{prop}
\begin{proof} The case of localization was already observed in Examples
\ref{ex.rfmods}. If $R \to S$ is \'etale then so is $R^e \to R^e \tensor_R S$
and $R^e \to S^e$. Thus, by \cite[Corollaire 17.3.4]{EGA4.4} the relative
Frobenius is \'etale too. It is an isomorphism if it is an isomorphism on
fibers, thus we may assume $R=k$ is a field and $S$ is a finite product of
separable algebraic extensions of $k$. Then the claim is easy, see
\cite[Theorem A.1.4]{Eisenbud.CommAlg}

Let $R \to \wh{R}$ be the $I$-adic completion as in (2). By assumption, $R^e$
is a finitely generated right $R$--module. Therefore $R^e \tensor_R \wh{R} =
\invlim R^e/R^eI^t$ by \cite[Theorem 7.2]{Eisenbud.CommAlg}. Using that the
sequence $R^eI^t=I^{t[p^e]}$ is cofinal within the powers $I^t$ of $I$ we
conclude
\[
    R^e \tensor_R \wh{R} \cong \invlim\frac{R^e}{R^eI^t} =
    \invlim\frac{R^e}{I^{t[p^e]}R^e} \cong \wh{R}^e.
\]
\end{proof}
Generally, the property of being finitely generated is not preserved by
restriction. The following is an important exception:

\begin{prop}\label{prop.finitetale} Let $S$ be finite \'etale over $R_x$
with $x \in R$. Then, a finitely generated unit $S[F^e]$--module $\Mm$ is
finitely generated as a unit $R[F^e]$--module.
\end{prop}
\begin{proof}
 The module
finiteness of $S$ over $R_x$ together with Proposition \ref{prop.loccompet}
shows that $\Mm$ is finitely generated as a unit $R_x[F^e]$--module. Now
\cite[Proposition 6.8.1.]{Em.Kis} finishes the proof.
\end{proof}

\begin{cor}\label{cor.rootComp}
    Let $R$ be regular and $R \to S$ be one of the cases of Proposition \ref{prop.loccompet}, that is in particular, $R^e \tensor S \cong S^e$.
    Let $\Mm$ be a finitely generated unit $R[F^e]$--module with root $M \subseteq \Mm$. Let
    $\Nn$ be a finitely generated unit $S[F^e]$--submodule of $S \tensor_R
    \Mm$. Then $\Nn \cap M$ is a root of the finitely generated unit
    $R[F^e]$--module $\Nn \cap \Mm$.
\end{cor}
\begin{proof}
By assumption, the finitely generated unit $S[F^e]$--modules $S \tensor_R \Mm$
and $\Nn$ are unit $R[F^e]$--modules (though quite likely \emph{not} finitely
generated as $R[F^e]$--modules). The intersection of the two unit
$R[F^e]$--submodules $\Mm$ and $\Nn$ is a unit $R[F^e]$--module. As it is a
submodule of the finitely generated module $\Mm$, it follows that $\Mm \cap
\Nn$ is a finitely generated unit $R[F^e]$--module since the category of
finitely generated unit $R[F^e]$--modules is abelian.

To check that the finitely generated module $N \defeq \Nn \cap M$ is a root of
$\Nn$ means that $\bigcup F_{S}^{er}(N) = \Nn$ and $N \subseteq F_{S}^{e}(N)$.
Thus
\[
    F^e_R(\Nn \cap M) = F^e_R(\Nn) \cap F^e_R(M) = F^e_{S}(\Nn) \cap
    M \supseteq \Nn \cap M
\]
and
\[
    \bigcup F_R^{er}(\Nn \cap M) = \bigcup (F^{er}_{S}(\Nn) \cap M) =
    \Nn \cap M.
\]
The key point was that for the $S$--submodule $\Nn$ of $S \tensor
\Mm$ one has
\[
    F^e_R(\Nn) = F^e_{S}(\Nn)
\]
by assumption.
\end{proof}
It is important to keep in mind that we did not exclude the case that $\Mm \cap
\Nn$ is zero in the last corollary. In particular it follows that $\Nn \cap \Mm
= 0$ if and only if $\Nn \cap M = 0$. Also note that $N \defeq S \tensor M \cap
\Nn$ is a root of the $S[F^e]$--module $\Nn$ and naturally $N \cap M= \Nn \cap
M$ is a root of the $R[F^e]$--module $\Nn \cap \Mm$.

\subsection{The minimal root}\label{sec.minroot}

Building on the last proposition and corollary we prove a result on the
existence of minimal roots for regular, $F$--finite local rings. This was
previously only known in the complete case \cite[Theorem 3.5]{Lyub}. First we
recall two easily verifiable facts, namely that the intersection of finitely
many roots of a finitely generated unit $R[F^e]$--module $\Mm$ is a root of
$\Mm$, and that if $\Nn$ is a unit $R[F^e]$--submodule of $\Mm$, and $M$ a root
of $\Mm$, then $M \cap \Nn$ is a root of $\Nn$.
\begin{theorem}\label{thm.minroot}
    Let $R$ be regular local and $F$--finite and let $\Mm$ be a finitely
    generated unit $R[F^e]$--module. Then $\Mm$ has a unique minimal root.
\end{theorem}
\begin{proof}
Let $\wh{R}$ denote the $m$--adic completion. By \cite[Theorem 3.5]{Lyub} or
\cite[Proposition 2.20]{Manuel.PhD}, the finitely generated unit
$\wh{R}[F^e]$--module $\wh{R} \tensor \Mm$ has a unique minimal root $N$.
Proposition \ref{prop.loccompet} is also a unit $R[F^e]$--module. By Corollary
\ref{cor.rootComp} $M \defeq N \cap \Mm$ is a root of the unit $R[F^e]$--module
$\Mm \cap (\wh{R} \tensor \Mm) = \Mm$. Clearly, $\wh{R} \tensor M \subseteq N$.
Since $\wh{R} \tensor M$ is a root of $\wh{R} \tensor \Mm$, it contains $N$ by
minimality of $N$. Therefore $N = \wh{R} \tensor M$. Now it follows easily that
$M$ is indeed the unique minimal root of $\Mm$. If $M' \subseteq M$ is another
root, we have, by minimality of $N$ the inclusion of roots $N \subseteq \wh{R}
\tensor M' \subseteq \wh{R} \tensor M$ of $\wh{R} \tensor \Mm$. Since the first
and last are equal we have that $M'$ and $M$ are equal upon completion. Thus
$M=M'$ by faithfully flat descent.
\end{proof}
A consequence of the above proof is the following corollary.
\begin{cor}
  Let $(R,m)$ be regular local and $F$--finite. Let $\Mm$ be a finitely
  generated unit $R[F]$--module with minimal root $M$. Then $\wh{R}
  \tensor M$ is the minimal root of $\wh{R} \tensor \Mm$.
\end{cor}
The question of whether minimal roots exist for not necessarily
local rings $R$ remains open.

\subsection{Duality for $R[F^e]$--modules} A key tool in local algebra is
Matlis Duality. If $(R,m)$ is local then the Matlis dual functor is defined as
$D(\usc) \defeq \Hom(\usc, E_{R/m})$, where $E_{R/m}$ denotes the injective
hull of $R/m$. We seek to extend $D(=D_R)$ to a Functor from $R[F^e]$--modules
to $R[F^e]$--modules. How this can be done is described in \cite{Manuel.PhD},
Chapter 4, in complete detail, as a consequence of a general investigation on
how to extend contravariant functors to incorporate Frobenius action. Here we
only give the bare minimum to establish the extension of $D$, most of the
material can already be found in \cite{Lyub}, Section 4.

\begin{prop}\label{proplFcommutesD}
    Let\/ $R$ be regular complete and local. The natural map
    \[
        \psi_\Mm: \F[e]\Hom(\Mm,E_R) \to \Hom(\F[e]\Mm,\F[e]E_R)
    \]
    is an isomorphism if $R$ is $F$--finite or if $\Mm$ is finitely presented or
    cofinite.

    In these cases we have an isomorphism of functors $\Psi: D \circ \F[e] \cong \F[e] \circ
    D$, that is, Matlis duality commutes with Frobenius.
\end{prop}
\begin{proof} In the first two cases this follows from the fact that the
natural map
\[
    \psi: S \tensor \Hom_R(\Mm,\Nn) \to \Hom_S(S \tensor \Mm, S \tensor \Nn)
\]
is an isomorphism provided that $R \to S$ is a flat map or rings and either $S$
is module finite over $R$ \cite[Proposition 4.9]{Manuel.PhD} or $\Mm$ is
finitely presented \cite[Proposition 2.10]{Eisenbud.CommAlg}. Thus it remains
to treat the case that $\Mm$ is a cofinite $R$--module which is treated in
\cite[Lemma 4.1]{Lyub}.

Fixing a unit structure on $E_R$ by fixing an isomorphism with $H^n_m(R)$ and
combining the above isomorphism with this unit $R[F^e]$--structure on $E_R$ we
get a natural (after the unit $R[F^e]$-structure on $E_{R/m}$ is fixed)
isomorphism
\[
    D(\F[e]\Mm) \cong \Hom(\F[e]\Mm,\F[e]E_R) \cong \F[e]\Hom(\Mm,E_R) =
    \F[e]D(\Mm)
\]
as desired.
\end{proof}
Now assume that $R$ is complete. Let $\F[e]\Mm \to[\theta^e] \Mm$ be an
$R[F^e]$--module which is finitely generated or cofinite as an $R$--module or
assume that $R$ is $F$--finite.

Applying the Matlis dual functor $D(\usc) = \Hom(\usc,E_R)$ to the structural
morphism of $\Mm$ and composing with the isomorphism of Proposition
\ref{proplFcommutesD} one obtains a map
\[
    \beta^e: D(\Mm) \to[D(\theta^e)] D(\F[e]\Mm) \to[\Psi_\Mm] \F[e](D(\Mm))
\]
whose second part is just the isomorphism $\Psi$ form the last Proposition.
\begin{defn}
    Let $R$ be complete and $(\Mm,\theta^e)$ an $R[F^e]$--module (finite\-ly
    generated or cofinite as an $R$--module, if $R$ is not $F$--finite).
    If $\beta^e \defeq \Psi_\Mm \circ D(\theta^e)$, then
    \[
        \Dd(\Mm) \defeq \dirlim(D(\Mm) \to[\beta^e] \F[e]D(\Mm)
        \to[\protect{\F[e]\beta^e}] \F[2e]D(\Mm) \to \ldots \quad)
    \]
    is the unit $R[F^e]$--module generated by $\beta^e$. On
    $R[F^e]$--modules which are cofinite as $R$--modules this defines an
    exact functor.
\end{defn}
The exactness claim is clear since Matlis duality and direct limits are exact
functors. If $\Mm$ is a unit $R[F^e]$--module then $\Dd(\Mm)=D(\Mm)$, since
$\beta^e$ is an isomorphism in this case. If $\Mm$ is cofinite as an
$R$--module then $D(\Mm)$ is a finitely generated $R$--module. Therefore
$\Dd(\Mm)$ is a finitely generated unit $R[F^e]$--module, since $D(\Mm)$, its
generator, is a finitely generated $R$--module. If in addition $\theta^e$ is
surjective, then $\beta^e$ is injective and therefore $D(\Mm)$ is a root of
$\Dd(\Mm)$.

\begin{notation}\label{not.F}
We introduce some notation from \cite{HaSp}. An element $m \in
\Mm$ of the $R[F^e]$--module $(\Mm,\theta^e)$ is called
$F$--nilpotent if $F^{re}(m)=0$ for some $r$. Then $\Mm$ is called
\emph{$F$--nilpotent} if $F^{er}(\Mm)=0$ for some $r \geq 0$. It
is possible that every element of $\Mm$ is $F$--nilpotent but
$\Mm$ itself is not, since $F$--nilpotency for $\Mm$ requires that
all $m \in \Mm$ are killed by the \emph{same} power of $F^e$. In
particular the sub $F[R^e]$--module consisting of all
$F$--nilpotent elements $\Fnil{\Mm}$ need not be nilpotent in
general. If $\theta^e$ is surjective, then $\Mm$ is called
\emph{F--full}. Note that $F$--fullness does not mean $F^e$ is
surjective but merely that the submodule
$F^e(\Mm)=\theta^e(\F[e]\Mm)$ is all of $\Mm$. Finally we say that
$\Mm$ is \emph{$F$--reduced} if $F^e$ acts injectively.
\end{notation}

The above notions are the same if we view $\Mm$ as an
$R[F^{er}]$--module for some $r \geq 0$. Therefore they are valid
without reference to a specific $e$.

We are lead to some functorial constructions for
$R[F^e]$--modules. The $R[F^e]$--submodule consisting of all
$F$--nilpotent elements of $\Mm$ we denote by $\Fnil{\Mm}=\{\,m
\in \Mm\,|\,F^{er}(m)=0 \text{ for some } r\,\}$. The quotient
$\Mm/\Fnil{\Mm}$ is the biggest $F$--reduced quotient, we denote
it by $\Fred{\Mm}$. The $R[F^e]$--submodule $\Fful{\Mm}=\bigcap
F^{er}(\Mm)$ is the largest $F$--full submodule. If $\Mm$ is a
cofinite $R$--module, then the decreasing chain of
$R[F^e]$--submodules $F^{er}(\Mm)$ stabilizes and we have
$\Fful{\Mm}=F^{er}(\Mm)$ for some $r
> 0$. One can check that the operations\/ $\Fful{(\usc)}$ and\/ $\Fred{(\usc)}$
mutually commute which makes the\/ $F$--full and  $F$--reduced subquotient\/
$\Ffred{\Mm}=\Fred{(\Fful{\Mm})}= \Fful{(\Fred{\Mm})}$ of an\/
$R[F^e]$--module\/ $\Mm$ $F$-reduced and $F$-full.

The following summary (see \cite[Section 4]{Lyub} for proofs) of the most
important properties of the functor $\Dd$ shows the significance of the just
introduced notions in our context.
\begin{prop}\label{prop.Dd}
    Let\/ $(R,m)$ be a regular, complete\/ $k$--algebra and let\/ $\Mm$ be a
    $R[F^e]$--module that is cofinite as an\/ $R$--module. Then
    \begin{enumerate}
    \item $\Dd(\Mm)=0$ if and only if\/ $\Mm$ is\/ $F$--nilpotent.\/ If $\Nn$ is
    also a cofinite $R[F^e]$--module, then $\Dd(\Mm) \cong
        \Dd(\Nn)$ if and only if\/ $\Ffred{\Mm} \cong \Ffred{\Nn}$.
    \item If\/ $\Mm$ is\/ $F$--full, then\/ $D(\Mm)$ is a root of\/ $\Dd(\Mm)$. If\/ $\Mm$ is
        also\/ $F$--reduced, then\/ $D(\Mm)$ is the unique minimal root.
    \item Every unit\/ $R[F^e]$--submodule $\Mm'$ of\/ $\Dd(\Mm)$ arises as\/ $\Dd(N)$
        for some $R[F^e]$--submodule of\/ $\Mm$.
    \item $\Dd$ is an isomorphism between the lattice of graded\/
        $R[F^e]$--modules quotients of\/ $\Mm$ (up to $\Ffred{(\usc)}$) and the lattice of
        unit\/ $R[F^e]$--sub\-modules of\/ $\Dd(\Mm)$.
    \end{enumerate}
\end{prop}
\label{p.simplRF}As a final remark we point out that if $M$ is a non-zero
simple $F$--full $R[F^e]$--module, then $\Dd(M)$ is non-zero and therefore a
simple unit $R[F^e]$--module. This follows since a simple $R[F^e]$--module is
$F$--full if and only if it is $F$--reduced and therefore by the last
Proposition $\Dd(M)$ is simple (and automatically nonzero). If $F^e$ had a
kernel it would be a nontrivial $R[F^e]$--submodule and thus if $M$ is simple
the kernel of $F^e$ must be all of $M$. Thus $F^e(M)=0$ which contradicts the
$F$--fullness since this exactly means that $F^e(M)=M$.

\subsection{The main example: $H^d_m(A)$}\label{sec.HdRelateHi}
Let $(R,m)$ be complete regular local ring of dimension $n$, let $I$ be an
ideal of height $c=n-d$. We denote the quotient $R/I$ by $A$.

The top local cohomology module $H^d_m(A)$ is an $A[F^e]$--module (\cf Examples
\ref{ex.rfmods}) and, by restriction, also an $R[F^e]$--module. As an
$R[F^e]$--module it is generally not unit, but at least the structural map
\[
    {R^e \tensor_R H^d_m(R/I)}  \to  {H^d_m(R/I)}
\]
is surjective. This is equivalent to the map induced by the
projection
\[
    R/I^{[p^e]} \to R/I
\]
under the identification of ${R^e \tensor_R H^d_m(R/I)}$ with $H^d_m(R/I^{[p^e]})$. Thus, by
definition, $\Dd(H^d_m(R/I))$ is the limit of
\begin{equation}
    D(H^d_m(R/I)) \to D(H^d_m(R/I^{[p^e]})) \to D(H^d_m(R/I^{[p^{2e}]})) \to
    \ldots
\end{equation}
Using local duality \cite[Theorem 3.5.8]{BrunsHerzog} for the
complete, regular and local ring $R$ this directed sequence is
isomorphic to the following
\begin{equation}\label{eqn.DdfromExt}
    \Ext_R^{c}(R/I,R) \to \Ext_R^{c}(R/I^{[p^e]},R) \to \Ext_R^{c}(R/I^{[p^{2e}]},R)
    \to \ldots
\end{equation}
where, again, the maps are the ones induced from the natural
projections. Since the Frobenius powers of an ideal are cofinite
within the normal powers, we get that the limit of this sequence
is just $H^{n-i}_I(R)$. This is because an alternative definition
of $H^{n-i}_I(R)$ is as the right derived functor of the functor
$\Gamma_I(M)=\dirlim \Hom(R/I^{t},R)$ of sections with support in
$\Spec R/I$.\footnote{See \cite[Theorem 3.5.6]{BrunsHerzog} for
the equivalence with our definition of local cohomology via \Czech
complexes.} The not very serious issue on whether the unit
$R[F^e]$--structure on $H^c_I(R)$ coming from the computation via
$\Ext$'s is the same as the one coming from the \Czech complex is
dealt with in \cite{Lyub}, Propositions 1.8 and 1.11. Summarizing
we get:
\begin{prop}\label{prop.DdofHmisHi}
    Let\/ $(R,m)$ be regular, local, complete and\/ $F$--finite. Let\/ $A=R/I$
    for some ideal\/ $I$ of\/ $R$ of height $c=n-d$. Then
    \[
        \Dd(H^d_m(R/I)) \cong H^{c}_I(R)
    \]
    as unit\/ $R[F^e]$--modules.
\end{prop}
By definition of $\Dd(\usc)$, a root for $\Dd(H^d_m(A))$ is given
by
\begin{equation}\label{eqn.DdfromExt2}
    \beta^e: \Ext_R^{c}(R/I,R) \to \Ext_R^{c}(R/I^{[p^e]},R) \to[\cong] R^e \tensor \Ext_R^{c}(R/I,R)
\end{equation}
where the first part is induced from the surjection $R/I^{[p^e]} \to R/I$, and
the second is the isomorphism coming from the natural transformation $\Psi: R^e
\tensor \Hom(\usc,R) \cong \Hom(R^e \tensor \usc, R)$, cf.\ Proposition
\ref{proplFcommutesD}. It is straightforward that this natural transformation
for $\Hom$ induces a natural transformation on its right derived functors, the
$\Ext$'s.

If we drop the assumption of completeness in the preceding discussion, and
just assume that $(R,m)$ is local we still have that $H^c_I(R)$ arises as
the direct limit of
\[
    \Ext^c_R(R/I,R) \to \Ext^c_R(R/I^{[p^e]},R) \to \Ext^c_R(R/I^{[p^{2e}]},R)
    \to \cdots
\]
with maps induced from the natural projections $R/I^{[p^{er}]} \to
R/I^{[p^{e(r-1)}]}$. Together with the natural transformation
identifying $\Ext^c_R(R/I^{[p^{er}]},R)$ with $R^{er} \tensor
\Ext^c_R(R/I,R)$ this shows that $\Ext^c_R(R/I,R)$ is a generator
for $H^c_I(R)$. Upon completion we get the generator
\[
    \Ext^c_{\wh{R}}(\wh{R}/I\wh{R},\wh{R}) \to[\beta^e] \wh{R}^e \tensor \Ext^c_{\wh{R}}(\wh{R}/I\wh{R},\wh{R})
\]
of $H^c_{I\wh{R}}(\wh{R}) \cong \wh{R} \tensor H^c_I(R)$ where we freely used
the identification $\wh{R} \tensor \Ext^c_R(R/I,R) \cong
\Ext^c_{\wh{R}}(\wh{R}/I\wh{R},\wh{R})$.

\section{Brief tight closure review}

Tight closure is a powerful tool in commutative algebra introduced
by Mel Hochster and Craig Huneke about fifteen years ago
\cite{HH88}. There is a strong connection between the
singularities arising in the minimal model program and
singularities obtained from tight closure theory
\cite{Smith.sing}. One of the most significant is the equivalence
of the notions of rational singularity and $F$--rational type
which was established by Smith \cite{Smith.rat} and Hara
\cite{Hara}, and independently by Mehta and Shrinivas
\cite{MehtaSr}. The notion of $F$--rationality arises naturally
from tight closure: the local ring $(A,m)$ is called $F$--rational
if all ideals $I$ generated by a full system of parameters are
tightly closed, \ie $I = I^*$. In this section we briefly review
the tight closure theory needed for our local construction of
$\Ll(Y,X)$ given below. For a more detailed introduction to this
beautiful subject we recommend \cite{Smith.IntroTight,Hune.tight}
and later the more technical original papers \cite{HH90,HH89}.

Let $A$ be a noetherian ring. We denote by $A^\circ$ the subset of
elements of $r$ that are not contained in any minimal prime of $A$. Let $N
\subseteq M$ be a submodule of $M$. We denote by $N^{[p^e]}$ the image of
$\F[e]N$ in $\F[e]M$. The tight closure $N^*_M$ of $N$ inside of $M$ is
defined as follows:
\begin{defn}
Let $A$ be noetherian and $N \subseteq M$. The tight closure $N^*_M$ (or
just $N^*$ if $M$ is clear from the context) consists of all elements $m
\in M$, such that there exists a $c \in A^\circ$, such that for all $e \gg
0$
\[
    c \tensor m \in N^{[p^e]}.
\]
Here $N^{[p^e]}$ denotes the image of $\F[e]N$ in $\F[e]M$ and $c \tensor
m$ is an element of $\F[e]M$.
\end{defn}
If $N=I$ is just an ideal of $A$, the definition is much more transparent.
In this case $r \in A$ is in $I^*$ if and only if there is $c \in A^\circ$
such that $cr^{p^e} \in I^{[p^e]}$ for all $e \gg 0$. A module is tightly
closed if $N^* = N$. We have that $N \subseteq N^*$ as one expects from a
decent closure operation. If $N$ is noetherian, then $N^* = (N^*)^*$.
There are two related closure operations which are important for us.
\begin{defn}
    Let $N \subseteq M$ be $A$--modules. The \emph{finitistic tight closure} of
    $N$ inside of $M$ consists of all elements $m \in (N \cap M_0)^*_{M_0}$
    for some finitely generated $M_0 \subseteq M$. It is denoted by
    $N^{*fg}_M$.

    The \emph{Frobenius closure} $N^F_M$ consists of all elements $m \in
    M$ such that $1 \tensor m \in N^{[p^e]}$ for some $e \geq 0$.
\end{defn}
We immediately see that $N^{*fg} \subseteq N^*$ and that equality
holds if $M$ is finitely generated. Clearly, $N^F \subseteq N^*$.
For the zero submodule of the top local cohomology module of an
excellent, local, equidimensional ring $A$, the finitistic tight
closure is equal to the tight closure, \ie
$0^{*fg}_{H^d_m(A)}=0^*_{H^d_m(A)}$ (see \cite[Proposition
3.1.1]{SmithDiss}). In general, it is a hard question to decide if
the tight closure equals the finitistic tight closure, and it is
related to aspects of the localization problem in tight closure
theory (\cf \cite{LyubSmith.Comm}).

As our focus lies on modules with Frobenius actions we study the above
closure operations in this case more closely. The following is an
important proposition which is proved in \cite{LyubSmith.Comm},
Proposition 4.2.
\begin{prop}
    Let\/ $A$ be noetherian and let\/ $(M,\theta^e)$ be an\/ $A[F^e]$--module. If\/ $N$ is
    a\/
    $A[F^e]$--submodule, then so are\/ $N^*_M$, $N^{*fg}_M$ and\/ $N^F_M$.
\end{prop}
This is checked by observing that $(N^*)^{[p^e]} \subseteq
(N^{[p^e]})^*$. Then apply $\theta^e$ and use the easily
verifiable fact that $\theta^e(\usc^*) \subseteq \theta^e(\usc)^*$
to see that
\[
    F^e(N^*)=\theta^e((N^*)^{[p^e]}) \subseteq \theta^e((N^{[p^e]})^*)
    \subseteq (F^e(N))^* \subseteq N^*
\]
which finishes the argument. From this we get as an immediate corollary
that the tight closure of the zero $A[F^e]$--submodule is a Frobenius
stable submodule of any $A[F^e]$--module.
\begin{cor}\label{cor.0*isRFsub}
    Let\/ $A$ be a ring and let\/ $(M,F^e)$ be an\/ $A[F^e]$--module.
    Then\/
    $0^{*fg}_M$, $0^*_M$ and\/ $0^F_M=\Fnil{M}$ are\/ $A[F^e]$--submodules of\/ $M$.
\end{cor}
\subsection{Test ideals and test modules}
The elements ``$c$'' occurring in the definition of tight closure play a
special role. Those amongst them, that work for all tight closure tests
for all submodules of all finitely generated $A$--modules are called the
test elements of $A$.
\begin{defn}
    An element $c \in A^\circ$ is called a \emph{test element} if for all
    submodules $N \subseteq M$, of every finitely generated $A$--module
    $M$, we have $cN^*_M \subseteq N$. A test
    element is called \emph{completely stable test element} if its image
    in the completion of every local ring of $A$ is a test element.
\end{defn}
It is shown in \cite[Proposition 8.33]{HH90}, that it is enough to range
over all ideals of $A$ in this definition, \ie $c$ is a test element if
and only if for all ideals $I$ and all $x \in I^*$ we have $cx^{p^e} \in
I^{[p^e]}$ for all $e \geq 0$. Thus, the test elements are those elements
$c$ occurring in the definition of tight closure which work for all tight
closure memberships of all submodules of all finitely generated
$A$--modules. A nontrivial key result is that in most cases, test elements
(and even completely stable test elements) exist:
\begin{prop}\label{prop.TestEl}
    Let\/ $A$ be reduced and of finite type over an excellent local ring.
    Then\/ $A$ has completely stable test elements. Specifically, any element\/ $c
    \in A^\circ$ such that\/ $A_c$ regular has a power which is a completely
    stable test element.
\end{prop}
The proof of this is quite technical and can be found in Chapter 6
of \cite{HH89}. Results in lesser generality (for example, when
$A$ is $F$--finite) are obtained fairly easily: for a good account
see \cite{Smith.IntroTight,Hune.tight}.

The ideal $\tau_A$ generated by all test elements is called the test
ideal. As remarked, $\tau_A = \bigcap (I:_A I^*)$ where the intersection
ranges over all ideals $I$ of $A$. This naturally leads one to consider
variants of the test ideal by restricting the class of ideals this
intersection ranges over. The \emph{parameter test ideal} of a local ring
$(A,m)$ is the ideal $\wt{\tau}_A = \bigcap (I:_A I^*)$ where the
intersection ranges over all ideals generated by a full system of
parameters. If $A$ is Cohen--Macaulay, it follows from the definition of
$H^d_m(A)$ as $\dirlim A/(\xd)^{[p^e]}$ that $\wt{\tau}_A =
\Ann_A(0^*_{H^d_m(A)})$ \cite[Proposition 4.1.4]{SmithDiss} where $\xd$ is
a system of parameters for the local ring $(A,m)$. If $A$ is only an
excellent domain, then $\wt{\tau}_A \subseteq \Ann_A(0^*_{H^d_m(A)})$.
Further generalizing, the \emph{parameter test module} is defined as
$\tau_{\omega_A} = \Ann_{\omega_A} 0^*_{H^d_m(A)} = \omega_A \cap
\Ann_{\omega_{\wh{A}}}0^*_{H^d_m(\wh{A})}$ where the action of $\omega_A$
on $H^d_m(A)$ is the one coming from the Matlis duality pairing $H^d_m(A)
\times \omega_{\wh{A}} \to E_A$. Of course we require here that $A$ has a
canonical module.
\begin{lem}\label{lem.tauNonzero}
    Let\/ $A$ be reduced, excellent, local and equidimensional with canonical
    module\/ $\omega_A$. If\/ $c$ is a parameter test element, then\/
    $c\omega_A \subseteq \tau_{\omega_A}$. In particular,\/ $\tau_{\omega_A}$ is
    nonzero.
\end{lem}
\begin{proof}
Let $c$ be a parameter test element. In particular, $c$ annihilates the
finitistic tight closure of zero in $H^d_m(A)$. Therefore, for every $\phi
\in \omega_A$ and $\eta \in 0^*_{H^d_m(A)}=0^{*fg}_{H^d_m(A)}$ we have
$c\phi \cdot \eta = \phi \cdot (c\eta) = \phi \cdot 0 = 0$ where
``$\cdot$'' represents the Matlis duality pairing. This shows that
$c\omega_A \subseteq \tau_{\omega_A}$. The hypotheses on $A$ ensure by
\cite[Remark 2.2(e)]{HH.CanMod} that the canonical module is faithful, \ie
$c\omega_A \neq 0$. Therefore the last part of the lemma follows from the
existence of test elements (Proposition \ref{prop.TestEl}), since a test
element is also a parameter test element.
\end{proof}

\subsection{$F$-rationality and local cohomology}
The tight closure of zero in the top local cohomology module
$H^d_m(A)$ of a local ring $(A,m)$ plays a role as the obstruction
to $F$--rationality of $A$. Its distinguishing property is that it
is the maximal proper $A[F^e]$--submodule of $H^d_m(A)$. Precisely
the following is the case:
\begin{theorem}\label{thm.MaxRFofH}
    Let\/ $(A,m)$ be reduced, excellent and analytically irreduci\-ble.
    Then, the tight closure of zero,\/ $0^*_{H^d_m(A)}$, in\/ $H^d_m(A)$
    is the unique maximal proper\/ $A[F^e]$--submodule of\/ $H^d_m(A)$.

    The quotient\/ $H^d_m(A)/0^*_{H^d_m(A)}$ is a nonzero
    simple\/ $F$--reduced and\/ $F$--full\/ $A[F^e]$--module.
\end{theorem}
\begin{proof}
The case $e=1$ of the first part was shown by Smith in
\cite{SmithDiss}, Theorem 3.1.4. The case $e \geq 1$ can be
obtained similarly, see \cite{Manuel.PhD}, Theorem 5.9. Because
$0^*_{H^d_m(A)}$ is the maximal proper $A[F^e]$--submodule,
$H^d_m(A)/0^*_{H^d_m(A)}$ is a simple $A[F^e]$--module quotient.
It remains to show that it is $F$--reduced (a simple
$A[F^e]$--module is $F$--full if and only if it is $F$--reduced).
For this note that the kernel of $F$ is a $A[F^e]$--submodule and,
by simplicity, it must either be zero ($F$--reduced) or all of
$H^d_m(A)/0^*_{H^d_m(A)}$. In the second case, this implies that
$F(H^d_m(A)) \subseteq 0^*_{H^d_m(A)}$. Since $H^d_m(A)$ is a unit
$A[F^e]$--module (enough that the structural map $\theta$ is
surjective) we have that $F(H^d_m(A)) = H^d_m(A)$. This
contradicts the fact that $0^*_{H^d_m(A)}$ is a proper submodule.
Thus the quotient is $F$--reduced and $F$--full.
\end{proof}
To avoid the assumption of analytically irreducible we give a
version of the above for the case that $A$ is an excellent
equidimensional ring. As the statement is about $H^d_m(A)$ which
does not discriminate between $A$ and its completion, we state the
result for a complete $A$; in general one has to take the minimal
primes of the completion of $A$ in the statement below.
\begin{cor}\label{cor.maxAFsubs}
    Let $A$ be a complete, local, reduced and equidimensional ring
    of dimension $d$. Let $P_1, \ldots, P_k$ be the minimal primes
    of $A$. Then the maximal proper $A[F^e]$--submodules are
    precisely
    \[
        M_i \defeq \ker( H^d_m(A) \to
        \frac{H^d_m(A/P_i)}{0^*_{H^d_m(A/P_i)}})
    \]
    where $i=1\ldots k$. Furthermore, the tight closure of zero, $0^*_{H^d_m(A)}$, in
    $H^d_m(A)$ is the intersection of all maximal proper
    $A[F^e]$--submodules. Even though $H^d_m(A)/0^*_{H^d_m(A)}$
    night not be simple as an $A[F^e]$--module, it is still $F$--full
    and $F$--reduced.
\end{cor}
\begin{proof}
Since tight closure can be checked modulo minimal primes, the last
statement is immediate.\footnote{This is generally proved for the
tight closure of ideals in the literature (see \cite{HH90},
Proposition 6.25), but the same proof can be adapted for
submodules.} By the last Theorem $0^*_{H^d_m(A/P_i)}$ is the
maximal proper $A[F^e]$--submodule. Thus $H^d_m(A)/M_i \cong
H^d_m(A/P_i)/0^*_{H^d_m(A/P_i)}$ is simple. Thus $M_i$ is a
maximal proper $A[F^e]$--submodule. To check that the $M_i$'s are
all the maximal proper $A[F^e]$--submodule let $M$ be a
$A[F^e]$--submodule of $H^d_m(A)$ not contained in any of the
$M_i$. This implies that for all $i$ the image of $M$ in
$H^d_m(A/P_i)=H^d_m(A)/P_iH^d_m(A)$ is all of $H^d_m(A/P_i)$ (it
is an $A[F^e]$--module not contained in $0^*_{H^d_m(A/P_i)}$, thus
must be all of $H^d_m(A/P_i)$ by last Theorem). But this implies,
by the following lemma, that $M=H^d_m(A)$ and we are done.

It remains to remark that a possible kernel of $F^e$ on
$H^d_m(A)/0^*_{H^d_m(A)}$ would reduce to all of
$H^c_m(A/P_i)/0^*_{H^c_m(A/P_i)}$ for some $i$ and therefore imply
that $F(H^c_m(A/P_i))\subseteq 0^*_{H^c_m(A/P_i)}$ which is a
contradiction to $F$--fullness of $H^c_m(A/P_i)$.
\end{proof}
\begin{lem}
    Let $A$ be a noetherian ring, and let $M \subseteq H$ be an
    $A$--module such that for every minimal prime $P$ of $A$ one
    has $M+PH=H$. Then $M = H$.
\end{lem}
\begin{proof}
One immediately reduces to the case $M=0$. Successive application
of the assumption $H=PH$ implies that
\[
    H=(P_1\cdot\ldots\cdot P_k)^nH
\]
where the $P_i$'s are the minimal primes. But for large enough
$n$, a power of the product of all minimal primes is zero, thus
$H=0$.
\end{proof}

If $A$ is Cohen--Macaulay, the vanishing of the tight closure of
zero in $H^d_m(A)$ characterizes $F$--rationality of $A$ by
\cite{Smith.rat}, Theorem 2.6. By definition, $A$ is called
$F$--rational if and only if every ideal that is generated by a
system of parameters is tightly closed.

\section{The intersection homology module}
First we give a detailed proof of the main existence theorem as
sketched in the introduction.
\begin{theorem}\label{thm.main}
    Let $X$ be an irreducible smooth $k$--scheme, essentially of finite
    type over $k$, and let $Y$ be a closed irreducible subscheme of
    codimension $c$. Then $\Hh^c_Y(X,\Oo_X)$
    has a unique simple $\Dd_X$--submodule $\Ll(Y,X)$. This submodule is
    also the unique simple $\Oo_X[F^e]$--module and agrees with
    $\Hh^c_Y(X,\Oo_X)$ on the complement of any closed set containing
    the singular locus of\/ $Y$.
\end{theorem}
\begin{proof}
Write $Z=\op{Sing} Y$ and $Y' = Y-Z$ and $X'=X-Z$ and denote the
open inclusion $X' \subseteq X$ by $i$. First we assume that the
characteristic of $k$ is positive; at the end of the proof we
indicate how the proof is adapted to characteristic zero.

We first show that $H^c_{Y'}(X', \Oo_{X'})$ is simple as a unit
$\Oo_X[F^e]$--module: Quite generally we note that, $\Oo_{Y'}$ is a simple unit
$\Oo_{Y'}[F^e]$--module by observing that a nontrivial ideal $\Ii \subseteq
\Oo_{Y'}$ is never a unit submodule as the containment $\Ii^{[p^e]} \subseteq
\Ii$ is strict, \cf Examples \ref{ex.rfmods}. Using that $Y'$ is smooth and
irreducible, it follows that $\Oo_{Y'}$ is also simple as a $\Dd_{Y'}$--module.
This can be reduced, by \'etale invariance of $\Dd_Y'$, to the case
$Y'=\Spec(k[x_1,\ldots,x_d])$ where one can check it by hand. Under Kashiwara's
equivalence for $\Dd_X$--modules (\cite{Haastert.DirIm}), and for unit
$\Oo_X[F^e]$--modules (\cite{Em.Kis}, Theorem 5.10.1 or \cite{Lyub},
Proposition 3.1), the module $\Oo_{Y'}$ corresponds to
$\Hh^c_{Y'}(X',\Oo_{X'})$ (\cf \cite{Em.Kis}, Example 5.11.6). Therefore,
$\Hh^c_{Y'}(X',\Oo_{X'})$, is a simple $\Dd_{X'}$--module (simple unit
$\Oo_{X'}[F]$--module, respectively).

Since $\Hh^c_Y(X,\Oo_X)$ is a locally finitely generated unit
$\Oo_X[F^e]$--module and therefore, by \cite{Lyub}, Theorem 5.6,
it has finite length as a $\Dd_X$--module. This assures the
existence of simple $\Dd_X$--submodules of the $D$-module $\Hh^c_Y(X,\Oo_X)$ and
let $\Ll_1$ and $\Ll_2$ be two such. Observe the exact sequence
(see \cite{Ha.LocCohom}, Chapter 1)
\[
    0=\Hh^c_Z(X,\Oo_X) \to \Hh^c_Y(X,\Oo_X) \to
    \Hh^c_{Y'}(X,\Oo_X) \cong i_*\Hh^c_{Y'}(X',\Oo_{X'})
\]
where the last isomorphism is by excision and the vanishing of the
first module is because the codimension of $Z$ in $X$ is strictly
bigger than $c$. From this it follows that $\Hh^c_Y(X,\Oo_X)$ and
therefore $\Ll_i$ are submodules of $i_*\Hh^c_{Y'}(X',\Oo_{X'})$.
By adjointness of restriction and extension we have
\[
    0 \neq Hom_{\Oo_X}(\Ll_i,i_*\Hh^c_{Y'}(X',\Oo_{X'})) \cong
    \Hom_{\Oo_{X'}}(\Ll_i|_{X'},\Hh^c_{Y'}(X',\Oo_{X'}))
\]
which shows that $\Ll_i|_{X'}$ are nonzero submodules of
$\Hh^c_{Y'}(X',\Oo_{X'})$. By simplicity of the latter all three
have to be equal. In particular, the intersection of $\Ll_1$ with
$\Ll_2$ is nonzero. As both are simple, this implies that
$\Ll_1=\Ll_2=\Ll(Y,X)$ as claimed. Furthermore, since
$F^e(\Ll(Y,X))$ is also simple, it follows from the uniqueness
that $F^e(\Ll(Y,X))=\Ll(Y,X)$ and therefore $\Ll(Y,X)$ is also the
unique simple $\Oo_X[F^e]$--submodule for all $e$.

Essentially the same proof works in characteristic zero. The key
fact then is that $\Hh^c_Y(X,\Oo_X)$ is a holonomic
$\Dd_X$--module and that holonomic modules have finite length.
Also observe that for the smooth $Y'$ the structure sheaf
$\Oo_{Y'}$ is $\Dd_{Y'}$--simple which is well known and easy to
check by hand for the case $Y'=\Spec(k[x_1,\ldots,x_d])$. Then
Kashiwara's equivalence implies that the corresponding
$\Hh^c_{Y'}(X',\Oo_X')$ is a simple $\Dd_{X'}$--module. For all of
the above statements, see \cite{Borel.Dmod}.
\end{proof}
This proof is pretty much identical for zero and positive
characteristic. The metaresults which are used though are proved
by very different techniques in each case.

If $Y$ is not irreducible, then the above result breaks down since then
$\Hh^c_{Y'}(X',\Oo_{X'})$ is no longer simple. In the case that
$Y$ is equidimensional one can give a complete description of the
simple submodules of the module $\Hh^c_Y(X,\Oo_X)$. They correspond to the
irreducible components $Y_1,\ldots, Y_k$ of $Y$. For each
component we have an inclusion
\[
    \Ll(Y_i,X) \subseteq H^c_{Y_i}(X,\Oo_X) \subseteq
    H^c_Y(X,\Oo_X)
\]
which establishes $\Ll(Y_i,X)$ as simple submodules of
$H^c_Y(X,\Oo_X)$. That the right map is an inclusion uses
equidimensionality and follows from \cite{Ha.LocCohom},
Proposition 1.9 and Chapter 3. To see that these are all the
simple submodules of $\Hh^c(X,\Oo_X)$ we show that any submodule
$\Nn$ of $H^c_Y(X,\Oo_X)$ does contain one of the $\Ll(Y_i,X)$. At
least for one $i$, the restriction of $\Nn$ to a open subset of
$X$ containing $Y_i$ but none of the other components is a nonzero
submodule of $H^c_{Y_i}(X,\Oo_X)$ (using excision). But then $\Nn$
clearly contains $\Ll(Y_i,X)$ by its uniqueness and simplicity. We
get as a corollary:
\begin{cor}
    Let $Y$ be an equi-dimensional and reduced sub-scheme of the smooth
    $k$--variety $X$ of codimension $c$. Let $Y=Y_1 \cup \ldots \cup Y_k$
    be its decomposition into irreducible components. Then the
    simple $\Dd_X$--submodules of $H^c_Y(X,\Oo_X)$ are precisely
    the $\Ll(Y_i,X)$ for $i=0,\ldots,k$. In this case we denote by
    $\Ll(Y,X)$ the (direct) sum of all the $\Ll(Y_i,X)$.
    Furthermore, away from the singular locus of $Y$ we have that
    $\Ll(Y,X)$ agrees with $\Hh^c_Y(X,\Oo_X)$.
\end{cor}
The similarity of Theorem \ref{thm.MaxRFofH} and Corollary
\ref{cor.maxAFsubs} with Theorem \ref{thm.main} and the last
corollary was the original motivation which lead to the
construction of $\Ll(Y,X)$ which is given in the next section.

Note that by the uniqueness of $\Ll(Y,X)$ it is clear that it
localizes, \ie if $U$ is a open subset of $X$ then $\Ll(Y,X)|_U =
\Ll(Y\cap U, U)$.

In the case of positive characteristic, we state the following
slightly stronger local version of the last theorem.
\begin{theorem}\label{thm.mainlocal}
    Let $R$ be a regular local ring of positive characteristic.
    Let $A=R/I$ be equidimensional of codimension $c$ in $R$.
    Then the sum of all simple unit submodules of $H^c_I(R)$ is
    \[
        \Ll(A,R) = \oplus \Ll(R/P_i,R)
    \]
    where the $P_i$ range over the primes minimal over $I$ and
    $\Ll(R/P_i,R)$ is the unique simple unit
    $\Oo_X[F]$--submodule of $H^c_{P_i}(R)$. Moreover, if $f \in
    R$ is such that $A_f$ is regular, then $\Ll(A,R)_f \cong
    H^c_I(R)_f$. If $R$ is $F$--finite the same holds for the
    $D_R$--module structure.
\end{theorem}
\begin{proof}
The local results of \cite{Lyub} are not restricted to finitely
generated algebras over a field. With the assumption above
(sufficient to ensure that a finitely generated unit
$R[F^e]$--module has finite length as such) the previous proof
goes through verbatim.
\end{proof}

\subsection{Local construction of $\Ll(A,R)$ in positive characteristic.}
\label{sec.LocConstLl} From now on we assume that $k$ is of
positive characteristic and $R$ is $F$--finite. With the last
theorem we can dispose of the $D_R$--structure and entirely work
with the $R[F^e]$--structure in our investigation $\Ll(A,R)$. As
the construction of $\Ll(Y,X)$ which we are about to present is
local in nature, the language is adjusted accordingly. Moreover,
the local construction of $\Ll(A,R)$ can be reduced to the
complete case with help of the results on minimal roots of Section
\ref{sec.minroot}. Thus we assume for now that $(R,m)$ is a
complete, regular local and $F$--finite, and that $A=R/I$
equidimensional of codimension $c$. The philosophy behind the
description of the simple unit $R[F^e]$--module $\Ll(A,R)$ is to
identify its minimal root. As it turns out, the minimal root of
$\Ll(A,R)$ is the parameter test module, which, under Matlis
duality, corresponds by definition to the tight closure of zero
$0^*_{H^d_m(A)}$ in $H^d_m(A)$.
\begin{theorem}\label{thm.copLl}
    Let $(R,m)$ be a complete regular local and $F$--finite
    and let $A=R/I$ be equidimensional and
    of codimension $c$. Then, we have that
    \[
        \Ll(A,R) = \Dd(H^d_m(A)/0^*_{H^d_m(A)})
    \]
    where $0^*_{H^d_m(A)}$ is the tight closure of zero in $H^d_m(A)$.
\end{theorem}
\begin{proof}
First assume that $A$ is a domain. By Theorem \ref{thm.MaxRFofH}
the unique simple $R[F^e]$--module quotient of $H^d_m(A)$ is
$H^d_m(A)/0^*_{H^d_m(A)}$. Therefore
\[
    \Dd(H^d_m(A)/0^*_{H^d_m(A)}) \subseteq \Dd(H^d_m(A)) \cong
    H^c_I(R)
\]
is a nonzero simple unit $R[F^e]$--submodule of $H^d_m(A)$. As
$\Ll(A,R)$ is the unique such they are equal.

If $A$ is only equidimensional, let $P_1,\ldots,P_k$ be its
minimal primes. In Corollary \ref{cor.maxAFsubs} we show that
\[
    0^*_{H^d_m(A)} = \ker( H^d_m(A) \to \oplus_{i=1}^k
    H^d_m(A/P_i)/0^*_{H^d_m(A/P_i)})
\]
Applying the functor $\Dd$ and using the domain case for $A/P_i$
as just proved one checks that
\[
    \Dd(H^d_m(A)/0^*_{H^d_m(A)}) = \oplus_{i=1}^k \Ll(A/P_i,R) = \Ll(A,R)
\]
where the last equality is by definition.
\end{proof}

A more careful investigation of the construction of $\Ll(A,R)$ via
the duality functor $\Dd$ shows its connection with the parameter
test module. By definition, $\Dd(H^d_m(A)/0^*)=\Ll(A,R)$ is the
unit $R[F^e]$--module generated by the Matlis dual of the
$R[F^e]$--module structure on $H^d_m(A)/0^*_{H^d_m(A)}$. The
Matlis dual of $H^d_m(A)$ is the canonical module
$\omega_A=\Ext^c_R(A,R)$ of $A$. The dual of
$H^d_m(A)/0^*_{H^d_m(A)}$ is found as the annihilator of
$0^*_{H^d_m(A)}$ under the Matlis duality pairing $\omega_A \times
H^d_m(A) \to E_{R/m}$. By definition, this annihilator
$\Ann_{\omega_A} 0^*_{H^d_m(R)}$ is the parameter test module
$\tau_{\omega_A}$. Thus the inclusion of unit $R[F^e]$--modules
$\Ll(A,R) \subseteq H^c_I(R)$ arises as the limit of the following
map between their generators:
\begin{equation}
\begin{split}
\xymatrix{
    {\omega_A}\ar[r] & {R^e \tensor \omega_A}\ar[r] & {\cdots} & {} & { = H^c_I(R)} \\
    {\tau_{\omega_A}}\ar@{^(->}[r]\ar@{^(->}[u] & {R^e \tensor \tau_{\omega_A}}\ar@{^(->}[r]\ar@{^(->}[u] & {\cdots} & {} & { = L(A,R)}\ar@{^(->}[u] \\
}
\end{split}
\end{equation}
Since, $H^d_m(A)/0^*_{H^d_m(A)}$ is $F$--reduced and $F$--full,
the bottom map $\tau_{\omega_A} \into R^e \tensor \tau_{\omega_A}$
is, by Proposition \ref{prop.Dd}(2), the unique minimal root of
$\Ll(A,R)$.
\begin{prop}
    Let $(R,m)$ be complete regular local and $F$--finite. Let $A=R/I$ be
    equidimensional of codimension $c$ in $R$. Then the parameter test module
    $\tau_{\omega_A}$ is the unique minimal root of $\Ll(A,R)$.
\end{prop}
Now we drop the assumption that $R$ be complete and only assume it
to be regular, local and $F$--finite.
\begin{theorem}
    Let $(R,m)$ be regular local and $F$--finite. Let $A=R/I$ be
    a domain of codimension $c$. Then $\Ll(A,R)$, the unique
    simple unit $R[F^e]$--submodule of $H^c_I(R)$, has the parameter
    test module $\tau_{\omega_A}$ as its minimal root.

    Furthermore, upon completion, $\wh{R} \tensor \Ll(A,R) =
    \Ll(\wh{A},\wh{R})$ .
\end{theorem}

\begin{proof}
By Corollary \ref{cor.rootComp}, $\Nn \defeq L(\wh{A},\wh{R}) \cap
H^c_I(R)$ is a unit $R[F^e]$-submodule of $H^c_I(R)$. A root of
$\Nn$ is found by intersecting $\Ll(\wh{A},\wh{R})$ with the root
$\omega_A$ of $H^c_I(R)$. Again by Corollary \ref{cor.rootComp}
this intersection is equal to $\tau_{\omega_A} = \omega_A \cap
\tau_{\omega_{\wh{A}}}$ which is the nonzero parameter test module
by Lemma \ref{lem.tauNonzero}. By simplicity of $\Ll(A,R)$ it is
therefore contained in $L(\wh{A},\wh{R}) \cap H^c_I(R)$, or put
differently,
\begin{equation}\label{eqn.loc1}
    \wh{R} \tensor \Ll(A,R) \subseteq \Ll(\wh{A},\wh{R}).
\end{equation}
If $R$ were a domain we would be done. Unfortunately, the
completion of a domain might not be a domain, but at least it is
equidimensional. If $P_1,\ldots,P_k$ are the primes minimal over
$I\wh{R}$, then Theorem \ref{thm.mainlocal} states that
$\Ll(\wh{A},\wh{R})=\oplus \Ll(\wh{R}/P_i,\wh{R})$. It remains to
show that $\Ll(\wh{R}/P_i,\wh{R}) \subseteq \wh{R} \tensor
\Ll(A,R)$. This is, by $\Ll(\wh{R}/P_i,\wh{R})$ being the unique
simple submodule of $\Hh^c_{P_i}(\wh{R})$, equivalent to $\wh{R}
\otimes \Ll(A,R) \cap \Hh^c_{P_i}(\wh{R}) \neq 0$. To see this let
$f \in R-{P_i}$ such that $A_f$ is regular, then, by the last part
of Theorem \ref{thm.mainlocal},
\[
    \wh{R} \tensor \Ll(A,R)_f = \wh{R} \tensor H^c_I(R)_f
    = H^c_{I\wh{R}}(\wh{R})_f \supseteq H^c_{P_i}(\wh{R})_f \neq
    0
\]
which shows the reverse inclusion of (\ref{eqn.loc1}). The
statement about the parameter test module now follows from the
beginning of the proof as we just showed that $\Ll(A,R)=\Nn$.
\end{proof}
As done before this proof can be adjusted to work for
equidimensional $A$ from the start. For simplicity we treated the
domain case and will do so from now on.

\begin{remk}
This construction of $\Ll(A,R)$ from its minimal root
$\tau_{\omega_A}$ enables one to explicitly construct
$D_R$--module generators for $\Ll(A,R)$: The image of any element
of $\tau_{\omega_A}$ in $H^c_I(R)$ is a generator of $\Ll(A,R)$,
in particular, if $c \in R$ is a test element such that $c^n\cdot
\eta \neq 0$ for $\eta \in \omega_A$, then $c \cdot \eta$
generates $\Ll(A,R)$ as a $D_R$--module.
\end{remk}

As another consequence of Corollary \ref{cor.rootComp} one sees
that the minimal root of $\Ll(A,R)$ is $\tau_{\omega_{\wh{A}}}$,
the minimal root of $\Ll(\wh{A},\wh{R})$, intersected with
$\omega_A$, the root of $H^c_I(R)$. By definition, this is the
parameter test module $\tau_{\omega_A}$ of $A$. Using Theorem
\ref{thm.minroot} it follows immediately that the parameter test
module commutes with completion, which was to the best of our
knowledge, unknown until now. We state this as a Proposition.
\begin{prop}
    Let $A$ be a domain which is a quotient of a regular local and
    $F$--finite ring. Then the parameter test module commutes
    with completion, \ie    $\tau_{\omega_{\wh{A}}} =
    \wh{A} \tensor \tau_{\omega_A}$.
\end{prop}

\subsection{Simplicity criteria for $H^c_I(R)$}
\label{sec.SimCritH}

With the connection between the tight closure of zero in
$H^d_m(A)$ and the simple unit $R[F^e]$--module $\Ll(A,R)$ just
derived, the characterization of Smith showing that
$0^*_{H^d_m(A)}$ governs the $F$--rationality of $A$, easily
implies a simplicity criteria for $H^c_I(R)$.
\begin{theorem}
    Let\/ $R$ be regular local and\/ $F$--finite. Let\/ $I$ be an ideal such
    that\/ $A=R/I$ is a domain. Then\/ $H^c_I(R)$ is\/
    $D_R$--simple if and only if the tight closure of zero in\/ $H^d_m(A)$ is\/
    $F$--nilpotent.
\end{theorem}
\begin{proof}
$H^c_I(R)$ is $D_R$--simple if and only if it is equal to $\Ll(A,R)$. Then
$\wh{R} \tensor \Ll(A,R)=\Dd(H^d_m(A)/0^*_{H^d_m(A)})$ is all of $\wh{R}
\tensor H^c_I(R)$ if and only if $\Dd(0^*_{H^d_m(A)})=0$, by exactness of
$\Dd$. This is the case if and only if $0^*_{H^d_m(A)}$ is $F$--nilpotent by
Proposition \ref{prop.Dd}.
\end{proof}
\begin{cor}
    Let\/ $R$ be regular, local and\/ $F$--finite. Let $A=R/I$ be a domain
    of codimension $c$. If\/ $A$ is\/ $F$--rational,
    then\/ $H^c_I(R)$ is\/  $D_R$--simple. If\/ $A$ is\/ $F$--injective
    (\ie $F$ acts injectively on\/ $H^d_m(A)$), then\/ $A$ is\/ $F$--rational
    if and only if\/ $H^c_I(R)$ is\/
    $D_R$--simple.
\end{cor}
\begin{proof}
By \cite{Smith.rat}, Theorem 2.6,  $F$--rationality of $A$ is
equivalent to $0^*_{H^d_m(A)} = 0$. Therefore, by the last
theorem, $\Ll(A,R)=H^c_I(R)$ if $A$ is $F$--rational. Conversely,
if $\Ll(A,R)=H^c_I(R)$ then $0^*_{H^d_m(A)}$ is $F$--nilpotent.
Under the assumption the $H^d_m(A)$ is $F$--reduced this implies
that $0^*_{H^d_m(A)}=0$, therefore $A$ is $F$--rational.
\end{proof}

This should be compared to the following characterization of
$F$--re\-gu\-la\-ri\-ty in terms of $D_A$--simplicity due to Smith:
\begin{prop}[\protect{\cite[2.2(4)]{SmithDmod}}]
    Let\/ $A$ be an\/ $F$--finite domain which is\/ $F$--split. Then\/ $A$ is
    strongly\/ $F$--regular if and only if\/ $A$ is simple as a\/
    $D_A$--module.
\end{prop}
Note that this proposition is a statement about the $D_A$--module structure of
$A$, \ie a statement about the differential operators on $A$ itself. This is
different from our approach as we work with the differential operators $D_R$ of
the regular $R$. Nevertheless, the similarity of the result is striking and
should be understood from the point of view Kashiwara's equivalence, \ie the
$D_A$--module $A$ should be studied via the corresponding $D_R$--module
$H^c_I(R)$.

We reformulate the simplicity criterion of $H^c_I(R)$ such that it
is a criterion solely on $A$, not referring to $H^d_m(A)$.
\begin{theorem}\label{prop.HcIsimpleCrit}
    Let\/ $R$ be regular, local and\/ $F$--finite. Let\/ $I$ be an ideal such
    that\/
    $A=R/I$ is a domain. If for all parameter ideals of\/ $A$ we
    have\/
    $J^F=J^*$, then\/ $H^c_I(R)$ is\/ $D_R$--simple.

    If\/ $A$ is Cohen--Macaulay, then\/ $H^c_I(R)$ is\/ $D_R$--simple if and only
    if\/
    $J^*=J^F$ for all parameter ideals\/ $J$.
\end{theorem}
\begin{proof}
We show that if $J^*=J^F$ for all parameter ideals, then $0^*_{H^d_m(A)}$ is
$F$--nilpotent, \ie $0^*_{H^d_m(A)}=0^F_{H^d_m(A)}$. Let $\eta \in H^d_m(A)$ be
represented by $z + (\xd)$ for some parameter ideal $J = (\xd)$, thinking of
$H^d_m(A)$ as the limit $\dirlim A/J^{[p^e]}$. Then the colon capturing
property of tight closure shows that $z \in 0^*_{H^d_m(A)}$ if and only if $z
\in J^*$ (\cf \cite[Proposition 3.1.1 ]{SmithDiss}). By our assumption
$J^*=J^F$, this implies that $z^{p^e} \in J^{[p^e]}$ for some $e > 0$.
Consequently, $F^e(\eta)=z^{p^e}+J^{[p^e]}$ is zero and thus every element of
$0^*_{H^d_m(A)}$ is $F$--nilpotent.

Under the assumption that $A$ is Cohen--Macaulay the same argument
can be reversed using that the limit system defining $H^c_m(A)$ is
injective.
\end{proof}

As we were dealing with the domain case above we remark that in
most cases when $A$ is not analytically irreducible the equality
$\Ll(A,R) = H^c_I(R)$ cannot hold, in particular $H^c_I(R)$ cannot
be simple.
\begin{prop}\label{prop.analrednotsimp}
Let $A=R/I$ be equidimensional, local and satisfy Serre's $S_2$
condition. Suppose that $A$ is \emph{not} analytically
irreducible, then $\Ll(A,R) \neq H^c_Y(A)$.
\end{prop}
\begin{proof}
    Equidimensionality and $S_2$--ness implies by
    \cite{HH.CanMod}, Corollary 3.7, that $H^d_m(A)$ is
    indecomposable. The properties of the duality functor $\Dd$
    show now that $\Dd(H^d_m(A))=H^c_{I\wh{R}}(\wh{R})$ is
    indecomposable as a unit $R[F^e]$--module. If $\wh{A}$
    is not a domain it follows, essentially by definition, that
    $\Ll(\wh{A},\wh{R})$ is decomposable, thus it cannot be equal
    to $H^c_{I\wh{R}}(\wh{R})$. As $\Ll$ behaves well under
    completion, it follows that $\Ll(A,R) \neq H^c_I(R)$.
\end{proof}

As an application of Theorem \ref{prop.HcIsimpleCrit} we extend
the last proposition to a  characterization of the simplicity of
$H^c_I(R)$ for the class of all domains $A$ which have only an
isolated singularity and whose normalization is $F$--rational. In
particular this yields a characterization of the simplicity of
$H^c_I(R)$ for $A=R/I$ a one dimensional domain.
\begin{prop}\label{prop.SmithEx}
    Let\/ $R$ be regular local and\/ $F$--finite. Let\/ $A=R/I$ be a local $S_2$
    domain with isolated singularity such that the normalization
    $\bar{A}$ is $F$--rational.
    Then\/ $H^c_I(A)$ is\/ $D_R$--simple if and only if $H^c_I(R)$ is analytically
    irreducible.
\end{prop}
\begin{proof}
If $H^c_I(R)$ is not analytically irreducible then, by Proposition
\ref{prop.analrednotsimp}, $H^c_I(R)$ is not $D_R$--simple
(equiv.\ unit $R[F^e]$--simple).

Since $R$ is excellent, $H^c_I(A)$ is analytically irreducible if
and only if $\bar{A}$ is local by \cite{EGA4.2}, (7.8.31 (vii)).
Let $z \in J^*$ for a parameter ideal $J=(\xd[y])$ of $A$. Then,
since $\bar{A}$ is $F$--rational and the expansion $\bar{J}$ of
$J$ to $\bar{A}$ is also a parameter ideal, one concludes that $z
\in \bar{J}^* = \bar{J}$. Let, for some $a_i \in \bar{A}$,
\[
    z = a_1y_1+ \ldots +a_dy_d
\]
be an equation witnessing this ideal membership. As we observe in
Lemma \ref{lem.IntClos} below, for some big enough $e$, all
$a^{p^e}$ are in $A$. Therefore $z^{p^e} = a_1^{p^e}y_1^{p^e}+
\ldots +a_d^{p^e}y_d^{p^e}$, which shows that $z^{p^e} \in
J^{[p^e]}$ since all $a^{p^e} \in A$. Thus $J^* = J^F$ and Theorem
\ref{prop.HcIsimpleCrit} implies that $H^c_I(R)$ is $D_R$--simple.
\end{proof}

\begin{lem}\label{lem.IntClos}
    Let\/ $(A,m)$ be a local domain with at worst isolated singularities.
    Then the normalization $\bar{A}$ of $A$ is local
    if and only if for all $x \in \bar{A}$ some power of $x$ lies in $A$.
\end{lem}
\begin{proof}
If there were $M_1$ and $M_2$, maximal ideals of $\bar{A}$ lying
over $m$, then, by assumption, for some $t \gg 0$ we have $(M_1)^t
\subseteq m \subseteq M_2$. Since $M_2$ is prime it follows that
already $M_1 \subseteq M_2$. The reverse inclusion follows by
symmetry and thus $M_1 = M_2$, therefore $\bar{A}$ is local.

Conversely, if $\bar{A}$ is local with maximal ideal $M$ it
follows that $\sqrt{m\bar{A}}=M$. Therefore, if $x\in M$ it
follows that $x^{n_0} \in m\bar{A}$ for sufficiently big $n_0$. We
want to conclude that $x^n \in m$ and thus is in $A$ for big
enough $n$. For this note that, by the assumption of isolated
singularity, the conductor ideal $C= (A:_{\bar{A}}\bar{A})$ is $M$
primary, \ie sufficiently high powers of $x$ lie in $C$. Now, if
$x^{n_0} \in m\bar{A}$ and $x^n \in C$, then
$x^{n_0+n}=x^{n_0}x^n$ is in $m$ itself. Finally, a unit $u$ of
$\bar{A}$ can be written as $u=\frac{ux}{x}$ for $x$ a nonunit of
$\bar{A}$. Since both, $x$ and $ux$ are not units, sufficiently
big powers are in $A$. Thus also sufficiently big powers of $u$
will be in $A$.
\end{proof}
In the case that $R$ is one--dimensional this yields a finite
characteristic analog of results of S.P. Smith
\cite{SmithSP.IntHom} and Yekutieli \cite{Yeku.ResDif}.
\begin{cor}\label{thm.curves}
    Let $A=R/I$ be a one--dimensional local domain with $R$ regular and
    $F$--finite. Then $H^c_I(R)$ is $D_R$--simple
    (equiv. unit $R[F^e]$--simple) if and only if $A$ is
    unibranch.
\end{cor}
\begin{proof}
As remarked in the proof of Proposition \ref{prop.SmithEx}, $A$ is
unibranch if and only if $A$ is analytically irreducible. As one
dimensional domains have at worst isolated singular points and
since the normalization is regular (and thus $F$--rational),
Proposition \ref{prop.SmithEx} applies.
\end{proof}

This last result that for curves $\Ll(A,R)$ is described in the same way in
positive characteristic as it is in characteristic zero is somewhat misleading.
In higher dimensions one expects that $\Ll(A,R)$ behaves significantly
different depending on the characteristic. For example, consider the ideal
$I=(xy-zw) \subseteq R = k[x,y,z,w]$. Then $A=R/I$ is the coordinate ring of
the cone over $\PP^1 \times \PP^1$ with the only singular point being the
vertex. The localization of $A$ at the vertex is $F$--rational. Therefore, our
results above shows that $H^1_I(R)$ is simple as a $D_R$--module in finite
characteristic. Nevertheless, in characteristic zero the module $H^1_I(R)$ is
not $D_R$--simple since the Bernstein-Sato polynomial of $xy-zw$ is
$(s-1)(s-2)$ and therefore has an integral zero of less than $-1$. This shows
that the $D_R$--submodule generated by $(xy-zw)^{-1}\in H^1_I(R)$ does not
contain $(xy-zw)^{-2}$. Therefore $H^c_I(R)$ has a proper $D_R$--submodule and
is therefore not $D_R$--simple.

This is in accordance with the Riemann--Hilbert type
correspondences in either characteristic. For zero characteristic,
the classical Riemann--Hilbert correspondence relates holonomic
$\Dd_X$--modules to constructible $\CC$--vectorspaces by means of
a vast generalization of de Rham theory, \ie to an ultimatively
topological theory. In positive characteristic, on the other hand,
the Emerton-Kisin correspondence relates finitely generated unit
$\Oo_X[F^e]$--modules to constructible $\FF^{p^e}$--sheaves on
$X$, generalizing Artin-Schreyer theory, which ultimately is a
coherent theory. This is one reason why there is no surprise for
the failure of a complete analogy of the description of the
intersection homology module $\Ll(Y,X)$ in positive and zero
characteristic.

\newcommand{\etalchar}[1]{$^{#1}$}
\providecommand{\bysame}{\leavevmode\hbox to3em{\hrulefill}\thinspace}
\providecommand{\MR}{\relax\ifhmode\unskip\space\fi MR }
\providecommand{\MRhref}[2]{%
  \href{http://www.ams.org/mathscinet-getitem?mr=#1}{#2}
}
\providecommand{\href}[2]{#2}


\begin{thebibliography}{BGK{\etalchar{+}}87}

\bibitem[BGK{\etalchar{+}}87]{Borel.Dmod}
A.~Borel, P.-P. Grivel, B.~Kaup, A.~Haefliger, B.~Malgrange, and F.~Ehlers,
  \emph{Algebraic ${D}$-modules}, Academic Press Inc., Boston, MA, 1987.

\bibitem[BH98]{BrunsHerzog}
Winfried Bruns and J{\"u}rgen Herzog, \emph{Cohen-{M}acaulay rings}, Cambridge
  University Press, Cambridge, UK, 1998.

\bibitem[BK81]{Bry.Kash}
J.-L. Brylinski and M.~Kashiwara, \emph{Kazhdan--{L}usztig conjecture and
  holonomic systems}, Invent. Math. \textbf{64} (1981), no.~3, 387--410.

\bibitem[Bli01]{Manuel.PhD}
Manuel Blickle, \emph{The intersection homology {$D$}--module in finite
  characteristic}, Ph.D. thesis, University of Michigan, 2001.

\bibitem[Bli03]{Manuel.DRF}
\bysame, \emph{The {$D$}--module structure of {$R[F]$}--modules}, Trans. Am.
  Math. Soc. \textbf{355} (2003), no.~4, 1647--1668.

\bibitem[Eis95]{Eisenbud.CommAlg}
David Eisenbud, \emph{Commutative algebra}, Springer-Verlag, New York, 1995.

\bibitem[EK99]{Em.Kis}
Mathew Emerton and Mark Kisin, \emph{Riemann--{H}ilbert correspondence for unit
  {$\mathcal{F}$}-crystals. {I}}, to appear, 1999.

\bibitem[EK00]{Em.Kis2}
\bysame, \emph{Riemann--{H}ilbert correspondence for unit
  {$\mathcal{F}$}-crystals. {I}{I}}, in preperation, 2000.

\bibitem[EK03]{Em.Kis3}
\bysame, \emph{{An introduction to the Riemann-Hilbert correspondence for unit
  {$\mathcal{F}$}-crystals}}, to appear in proceedings of Dwork memorial
  conference, 2003.

\bibitem[Gro65]{EGA4.2}
A.~Grothendieck, \emph{\'{E}l\'ements de g\'eom\'etrie alg\'ebrique. {I}{V}.
  \'{E}tude locale des sch\'emas et des morphismes de sch\'emas. {I}{I}}, Inst.
  Hautes \'Etudes Sci. Publ. Math. (1965), no.~24, 231.

\bibitem[Gro67]{EGA4.4}
\bysame, \emph{\'{E}l\'ements de g\'eom\'etrie alg\'ebrique. {I}{V}. \'{E}tude
  locale des sch\'emas et des morphismes de sch\'emas {I}{V}}, Inst. Hautes
  \'Etudes Sci. Publ. Math. (1967), no.~32, 361.

\bibitem[Haa88]{Haastert.DirIm}
Burkhard Haastert, \emph{On direct and inverse images of $d$-modules in prime
  characteristic}, Manuscripta Math. \textbf{62} (1988), no.~3, 341--354.

\bibitem[Har67]{Ha.LocCohom}
Robin Hartshorne, \emph{Local cohomology}, A seminar given by A. Grothendieck,
  Harvard University, Fall, vol. 1961, Springer-Verlag, Berlin, 1967.

\bibitem[Har98]{Hara}
Nobou Hara, \emph{A characterization of rational singularities in terms of the
  injectivity of the {F}robenius maps.}, American Journal of Mathematics
  \textbf{120} (1998), 981--996.

\bibitem[HH88]{HH88}
Melvin Hochster and Craig Huneke, \emph{Tightly closed ideals}, Bull. Amer.
  Math. Soc. (N.S.) \textbf{18} (1988), no.~1, 45--48.

\bibitem[HH89]{HH89}
\bysame, \emph{Tight closure and strong ${F}$-regularity}, M\'em. Soc. Math.
  France (N.S.) (1989), no.~38, 119--133, Colloque en l'honneur de Pierre
  Samuel (Orsay, 1987).

\bibitem[HH90]{HH90}
\bysame, \emph{Tight closure, invariant theory and the {B}riancon-{S}koda
  theorem}, Journal of the American Mathematical Society \textbf{3} (1990),
  31--116.

\bibitem[HH94]{HH.CanMod}
\bysame, \emph{Indecomposable canonical modules and connectedness}, Commutative
  algebra: syzygies, multiplicities, and birational algebra (South Hadley, MA,
  1992), Amer. Math. Soc., Providence, RI, 1994, pp.~197--208.

\bibitem[HS77]{HaSp}
Robin Hartshorne and Robert Speiser, \emph{Local cohomological dimension in
  characteristic $p$}, Annals of Mathematics \textbf{105} (1977), 45--79.

\bibitem[Hun96]{Hune.tight}
Craig Huneke, \emph{Tight closure and its applications}, Published for the
  Conference Board of the Mathematical Sciences, Washington, DC, 1996, With an
  appendix by Melvin Hochster.

\bibitem[LS01]{LyubSmith.Comm}
Gennady Lyubeznik and Karen~E. Smith, \emph{On the commutation of the test
  ideal with localization and completion}, 2001.

\bibitem[Lyu97]{Lyub}
Gennady Lyubeznik, \emph{$\mathcal{F}$-modules: an application to local
  cohomology and {$D$}-modules in characteristic $p>0$}, Journal f{\"ur} reine
  und angewandte Mathematik \textbf{491} (1997), 65--130.

\bibitem[MS97]{MehtaSr}
V.B. Mehta and V.~Srinivas, \emph{A characterization of rational
  singularities}, Asian Jounal of Mathematics \textbf{1} (1997), no.~2,
  249--271.

\bibitem[Smi88]{SmithSP.IntHom}
S.~P. Smith, \emph{The simple {$D$}-module associated to the intersection
  homology complex for a class of plane curves}, J. Pure Appl. Algebra
  \textbf{50} (1988), no.~3, 287--294.

\bibitem[Smi93]{SmithDiss}
Karen~E. Smith, \emph{Tight closure of parameter ideals and $f$--rationality},
  Ph.D. thesis, University of Michigan, 1993.

\bibitem[Smi95]{SmithDmod}
\bysame, \emph{The ${D}$-module structure of ${F}$-split rings}, Math. Res.
  Lett. \textbf{2} (1995), no.~4, 377--386.

\bibitem[Smi97a]{Smith.rat}
\bysame, \emph{${F}$-rational rings have rational singularities}, Amer. J.
  Math. \textbf{119} (1997), no.~1, 159--180.

\bibitem[Smi97b]{Smith.sing}
\bysame, \emph{Vanishing, singularities and effective bounds via prime
  characteristic local algebra}, Algebraic geometry---Santa Cruz 1995, Amer.
  Math. Soc., Providence, RI, 1997, pp.~289--325.

\bibitem[Smi01]{Smith.IntroTight}
\bysame, \emph{An introduction to tight closure}, Geometric and combinatorial
  aspects of commutative algebra (Messina, 1999), Dekker, New York, 2001,
  pp.~353--377.

\bibitem[Vil85]{Vil}
K.~Vilonen, \emph{Intersection homology ${D}$-module on local complete
  intersections with isolated singularities}, Invent. Math. \textbf{81} (1985),
  no.~1, 107--114.

\bibitem[Yek98]{Yeku.ResDif}
Amnon Yekutieli, \emph{Residues and differential operators on schemes}, Duke
  (1998), no.~95, 305.

\end{thebibliography}
\end{document}